\documentclass[12pt]{amsart}
\usepackage{amscd,amssymb,hyperref}

\usepackage{pstricks}

\usepackage{amssymb,latexsym}
\usepackage{amsmath}

\usepackage{amsthm,latexsym}
\usepackage{eucal,latexsym}

\theoremstyle{plain}

\theoremstyle{remark}

\theoremstyle{definition}

\numberwithin{equation}{section}

\begin{document}
\baselineskip=26pt
 \title {Open and closed string field theory interpreted in classical algebraic topology}
\author{Dennis Sullivan}
\dedicatory{Dedicated to Graeme Segal on his 60th birthday}
 \maketitle
\begin{center}

\end{center}

 \textit{Abstract}: There is an interpretation of open
string field theory in algebraic topology. An interpretation of
closed string field theory can be deduced from this open string
theory to obtain as well the interpretation of open and closed
string field theory combined. The algebraic structures derived
from the first string interactions are related to algebraic models
discussed in work of (Atiyah-Segal), (Moore-Segal) and (Getzler
and Segal). For example the Corollary 1 of \S1 says that the
homology of the space of paths in any manifold beginning and
ending on any submanifold has the structure of an associative
dialgebra satisfying the module or Frobenius compatibility (see
appendix). Corollary 2 gives another structure.

\S1\textit{Open string states in $M$}: The open string theory
interpretation in topology includes a collection of linear
categories $\vartheta M$ one for each ambient space $M$. The
objects of $\vartheta M$ are smooth oriented submanifolds
$L_a,L_b,L_c,...$ of $M$.  The set of morphisms $\vartheta_{ab}$
between two objects $L_a$ and $L_b$ are graded chain complexes,
linearly generated by smooth oriented families of paths from $L_a$
to $L_b$.  An element in $\vartheta_{ab}$ is called an open string
state.  A path is a piecewise smooth map [0,1]$\rightarrow M$.

The first open string interactions are
\newline $i$)\textit{two endpoint restrictions}:
    $\vartheta_{ab}\overset{r}\rightarrow \vartheta_{a'b}$ and $\vartheta_{ab} \overset{r}\rightarrow \vartheta_{ab'}$ where $L_{a'}$ is a submanifold of $L_a$ and $L_{b'}$ is
    a submanifold of $L_b$. Degree $r=$ $-$cod of submanifold.
    \newline $ii$)joining or composition
    $\vartheta_{ab}\otimes \vartheta_{bc}\overset{\wedge}\rightarrow
    \vartheta_{ac}$, degree $\wedge=$ $-$dim $L_b$
    \newline $iii$) cutting or cocomposition
    $\vartheta_{ac}\overset{\vee}\rightarrow \vartheta_{ab} \otimes
    \vartheta_{bc}$, degree $\vee=$ $-$cod $L_b$ +1

    Namely,
    \newline $i$)(restriction) for an open string state in
    $\vartheta_{ab}$(ie. a chain in $\vartheta_{ab}$) one can
    intersect transversally in $L_a$ the chain of beginning points in $L_a$ with $L_{a'}$
    to obtain a chain in $\vartheta_{a'b}$.  The same idea works in $L_b$
    for the endpoints of paths to construct
    $\vartheta_{ab} \overset{r}\rightarrow \vartheta_{ab'}$.
    \newline $ii$)(joining) the transversal intersection in $L_b$ of the chain of
    endpoints for an open string state in $\vartheta_{ab}$ with
    the chain of beginning points for an open string state in
    $\vartheta_{bc}$is a chain labelling composible paths which
    after composing defines an open string state in $\vartheta_{ac}$, and the composition
    $\vartheta_{ab} \otimes \vartheta_{bc}\overset{\wedge}{\rightarrow} \vartheta_{ac}$.
    \newline$iii$)(cutting) Now it is required that
    $L_a,L_b,L_c,...$ have oriented normal bundles.  For example,
    this is true if the ambient space $M$ is a smooth manifold.  Then
    given an $L_b$ and any open string state in $\vartheta_{ac}$
    we may transversally intersect in $M$ the paths with $L_b$.  The intersection chain
    labels cuttings of the path at $L_b$ defining
    $\vartheta_{ac}\overset {\vee}{\rightarrow} \vartheta_{ab} \otimes
    \vartheta_{bc}$. (We use Eilenberg-Zilber.)
    \newline  The operation $\vee$ refers to cutting at any time along the
    path whenever it crosses $L_b$.  We can also consider the
    operation $\vee_t$ of cutting at a specific time $t\epsilon[0,1]$.  All these
    $\vee_t$
    are chain homotopic.  In fact $\vee$ is the chain homotopy between
     $\vee_0$  cutting at time zero and $\vee_1$ cutting at time one.

    \textit{Remark}: Actually the above operations are directly defined by
    the above descriptions only for states satisfying
    transversality conditions. To go from such a typical definition to a complete
    definition perturbations of the identity creating
    transversality must be introduced. The combinatorics of these
    perturbations fits neatly into Stasheff's strong homotopy
    formalism [S].  An elegant treatment can be read in Fukaya et al
    [1], for the classical case of intersecting chains in a
    manifold.

    \textit{Theorem}: For each ambient oriented smooth manifold $M$ there is an open string
    category whose objects are smooth submanifolds
    $L_a,L_b,L_c,..$ and whose morphisms are chains
    $\vartheta_{\alpha\beta}$ on paths between objects $L_\alpha$
    and $L_\beta$.  Only the objects $L_a$ which are compact (without boundary) have identity maps (which
    commute with the boundary operator).
     For transversal open string states in
    $\vartheta_{\alpha\beta}$,...  composition $\wedge$ is associative,
    cocomposition $\vee$ is coassociative, and the derivation compatibility
    holds between $\vee$ and $\wedge(x,y)=x \cdot y$, $\vee (x \cdot y) = x \cdot \vee
    y$ +
     $\vee x \cdot y$(see appendix). $\wedge$ and $\vee_t$ commute with
    $\partial$ but $[\vee,\partial]=\vee_1-\vee_0$.

    On the full space of open string states, associativity for $\wedge$ and coassociativity for $\vee_t$ hold
     up to strong homotopy in the
    sense of Stasheff.  There are conjecturally similar
    strong homotopy statements for coassociativity of $\vee$ and the
    derivation or infinitesimal bialgebra compatibility between $\wedge$ and
    $\vee$.(see appendix).

    \textit{Corollary 1}: For each object $L_a$ the homology of
    $\vartheta_{aa}$ is an associative algebra via the composition
    operation $\wedge$ (with identity if $L_a$ is compact without boundary).
     The operation $\vee_t$ is a coassociative
    coalgebra (which if non zero implies $L_a$ cannot be deformed off of itself). The $\wedge$, $\vee_t$ dialgebra
    satisfies the module or Frobenius compatibility
    (see appendix).

    \textit{Proof of corollary}: $i$) The algebra statement
    follows from a) $\wedge$ commutes with $\partial$ operator on
    open string states and so passes to homology b)homotopy
    associativity at the chain level implies associativity at the
    homology level.
   \newline ii) a)The fixed time cutting operation $\vee_t$ also
    commutes with the $\partial$ operator and passes to homology.
    b) because different times are chain homotopic we can choose
    them conveniently to prove the module or Frobenius compatibility.  To
    calculate $\vee_t (x \cdot y)$ we can choose $t$ in $x$'s
    time to see that we get $\vee_t (x) \cdot y$ or in $y$'s
    time to see that we get $x \cdot \vee_t (y)$. See the remark
    2) for the rest.

    \textit{Sketch proof of theorem}: 1) One sees the indicated
    identities hold for transversal chains by looking at the
    picture.  For example, when cutting a joining of paths, the
    cut can happen in the first part or the second part. This
    yields the derivation compatibility.
    \newline 2) The strong homotopy properties follow using i)
    manifolds are locally contractible
    \newline ii) transversality can be created in manifolds by
    arbitrarily small pertubations.

    \textit{Remarks}: 1) The coalgebra $\vee_t$ is chain
    homotopic to $\vee_0$ which may be written as a composition
    involving the restriction and the diagonal mapping.  Let
    $L_{a'}$ be the transversal intersection of $L_a$ with itself.
     Then $\vee_0$ is the composition of, \textit{first the restriction of
     the beginning point} to $L_{a'}$, next the inclusion of
     $\vartheta_{a'a}$ into $\vartheta_{aa}$, \textit{next the diagonal
     map on generating chains of $\vartheta_{aa}$}, next the cartesian product on chains of the
     beginning point operator(thought of as a constant path) with the
     identity and finally Eilenberg-Zilber.  A similar composition and statement hold for
     $\vee_1$.
     \newline 2) We can use remark 1) to define a new coalgebra
     structure on homology when $L_a$ is deformable off itself,
     say to $L_b$.  Then define $\vee: \vartheta_{aa}
     \rightarrow \vartheta_{ab} \otimes \vartheta_{ba}$ cutting at variable time
     and note
     that $\vee_0$ and $\vee_1$ are zero on the chain level.
     Thus $\vee$ commutes with $\partial$ and passes to
     homology. We use the obvious equivalences
     $\vartheta_{aa}\sim \vartheta_{ba}\sim \vartheta_{ab}$ to
     obtain:

     \textit{Corollary 2}:  If $L_a$ is deformable off of itself,
     the homology of open string states on $L_a$ has the structure
     of an associative dialgebra satisfying the derivation or
     infinitesimal bialgebra compatibility (see appendix).

    \textit{Examples}:
    $i$)(\textit{manifolds}) $L_a=M$ the ambient space.  Then
    $\vartheta_{aa}$ is equivalent to the ordinary chains on $M$
    since {paths in $M$} is homotopy equivalent to $M$.  Then the
    strong homotopy associativity algebra structure on $\vartheta_{aa}$ is equivalent
    to the intersection algebra of chains on $M$.  The operation $ \vee_\circ \sim \vee_t \sim
    \vee_1$ is chain equivalent to the diagonal mapping on
    chains.  One recovers the known fact that on passing to homology one obtains a
    graded commutative algebra structure
    $C \otimes C \overset{\wedge}\rightarrow C$ and a graded cocommutative
    coalgebra structure $C \overset{\vee}\rightarrow C \otimes C$
    satisfying the module or Frobenius compatibility $\vee
    (x \cdot y)=x \cdot(\vee y)=(\vee x)\cdot y$ where the
    notation refers to multiplication on the left and right
    factors of the tensor product respectively (see appendix).

    Note when $M$ is a closed oriented manifold $\wedge$ and $\vee$
    are related by the non degenerate intersection pairing,
    Poincare duality.

    $ii$)(\textit{based loop space})$M$ is any space and $L_a$ is a point
    in $M$.  Then $\vartheta_{aa}$ is the chains on the based loop
    space of $M$ and the algebra structure on $\vartheta_{aa}$ is the
    Pontryagin algebra of chains on the based loop
    space (the original setting of Stasheff's work). No transversality is needed here
    because all paths are composible. Here one has Hopf's
    celebrated compatibility with the diagonal map $\vee'$ on
    chains that $\vee'$ is a map of algebras.  The connection
    of the latter with the open string theory here is a mystery
    (but compare [2] and remark 1) above).

    If $M$ is a manifold of dim $M$ near $L_a$ and $L_a$ is a point, the cocomposition
    $\vee_t$ is defined but is zero in
    homology.  The operation $\vee$ can then be refined to a
    chain mapping and passes to homology (remark 2)).  $\vartheta_{aa}$
    obtains a coassociative coalgebra structure on homology of
    degree (-dim $M$) +1 satisfying the derivation or infinitesimal bialgebra
    compatibility (of the theorem) with the Pontryagin product.
    Here one is splitting a based loop where it passes again
    through a (nearby) base point.

    $iii$) (\textit{free loop space}) Let $M= L$ x $L$ and $L_a \subset M$ be
    the diagonal.  Then paths in $M$ beginning and ending on $L_a$
    is homeomorphic to the free loop space of $L$= Maps (circle,
    $L$).  Then the algebra structure on $\vartheta_{aa}$ is
    chain homotopic to the loop product of "String Topology"[2].
    This is a graded commutative algebra structure on the homology
    of the free loop space of the manifold $L$.  The degree is
    zero if we grade by the negative codimension ($k-$dim$M$).

    The product interacts with the circle action differential
    $\triangle$ of degree +1.  The deviation of $\triangle$ from being a
    derivation of the loop product $\triangle (x
    \cdot y)-(\triangle x)\cdot y-x \cdot(\triangle y)$ is a Lie
    bracket of degree +1 which is compatible via the Leibniz
    identity with the loop product (all on homology).  This Lie bracket is a
    geometric version [2] of Gerstenhaber's bracket in the
    (Hochshild) deformation
    complex of an associative algebra. For simply connected
    closed manifolds $L$ the Hochshild complex $\oplus_k$ Hom ($A^{\otimes k},A$) applied to the intersection
    algebra $A$ of chains on $L$ is a model of the free loop space of
    $L$ (Cohen-Jones, Tradler) which realizes the above comparison (Tradler).

    The Lie product on the free loop space of degree +1 is
    compatible via the connecting morphism $M$ between equivariant
    homology and ordinary homology with a Lie bracket on the equivariant free loop space
    homology [2].  The latter Lie bracket generalizes to all manifolds
    the Goldman bracket (related to the Poisson structure on flat bundles
    over a surface) on the
    vector space generated by conjugacy classes in the fundamental
    group of a surface [Goldman] (see closed strings \S 2 below).

    If the coalgebra part $\vee_t$ of the Frobenius dialgebra on homology of the free loop space of $L$ is non zero, then
    $L$ is a closed
    manifold with non-zero Euler characteristic.  Otherwise a homotopy class
    of non-zero vector fields on $L$ allows a refining of the operation $\vee$
     cutting at variable time to an operation commuting with
     $\partial$ and we obtain in this case an infinitesimal
    bialgebra structure (appendix) on the homology of the free loop
    space.

    \S2\textit{Closed string states in $M$ (now called $L$)}:
    For closed string states in $L$ we take the chains for the
    equivariant free loop space of $L$ relative to the circle
    action rotating the domain.  There are maps
    \newline ...$\overset {C}
    \rightarrow$ closed string states in $L$
    $\overset{M}\rightarrow$ open string states on the diagonal in
    $L$x$L$ $\overset {E}\rightarrow$ closed string states in
    $L$ $\overset{C}\rightarrow$... \newline leading to the long exact sequence relating ordinary homology
    and equivariant circle homology.  Here we are thinking of the
    free loop space of $L$ as paths in $L$ x $L$ beginning and
    ending on the diagonal.

    The connecting chain map $C$ has degree -2 and intersects with a representative of the 1st chern
    class of the line bundle associated to the $S^1$ action made free by crossing with a contractible
    space on which $S^1$ acts freely.
    The chain map $M$ has degree +1 and is associated to adding a mark to a closed string
    in all ways to get a circle of free loops.  The chain map $E$ has
    degree zero and is associated to forgetting the mark on a loop to get a closed string.  The composition $EM=0$ and
    the composition $ME$ is $\triangle$ the differential
    associated to the circle action.

    The string product on closed string states satisfying Jacobi (at the transversal chain level)
     may
    be defined by the formula $[\alpha, \beta]= E(M \alpha \wedge
    M \beta)$ where $\wedge$ is the open string product (the procedure in example 3 above only satisfies
    Jacobi up to a non trivial chain homotopy).  Other
    independent closed string operations $c_n$ can be defined by
    $c_n (\alpha_1,\alpha_2,...,\alpha_n$)=$E(M \alpha_1 \wedge
    M \alpha_2 \wedge...\wedge M \alpha_n)$(cf.
    [2] and [G]).  These all commute with the $\partial$
    operator and satisfy other identities transversally [2].

    The collision operators $c_n$ pass to the reduced equivariant
    complex or \textit{reduced closed string} states which is
    defined to be the equivariant chain complex for the $S^1$
    pair, (free loop space, constant loops).

    We can define a closed string cobracket
    $s_2$ by the formula $s_2(\alpha)=(E \otimes E)(\vee (M \alpha))$.
    In the reduced complex $s_2$ commutes with $\partial$ and
    passes to homology (but not so in the unreduced complex).

    \textit{Theorem}: The closed string bracket
    $c_2(\alpha,\beta)= E(M\alpha \wedge M \beta)$ where $x \wedge y=\wedge(x \otimes y)$ and the closed
    string cobracket $s_2(\alpha)=(E \otimes E)(\vee M \alpha)$ satisfy
    respectively jacobi, cojacobi, and Drinfeld compatibility (appendix).
    The term satisy means either on the level of integral homology, for
    transversal chains on the chain level, or conjecturally at
    the Stasheff level of strong homotopy.

    \textit{Proof}: These formulae in terms of open strings are
    reinterpretations as in [2] of the definitions given in "Closed string
    operators in topology leading to Lie bialgebras and higher
    string algebra" [3]. There the
    identities at the transversal chain level were considered.

    \textit{Corollary}: Homology of reduced closed string states
    forms a Lie bialgebra, [3].

    \textit{Remark}: Independent splitting operations $s_3, s_4,...$
    can be defined similarly by iterations of $\vee$,
    $s_n(\alpha)=E \otimes...\otimes E(...\vee \otimes 1 \cdot \vee(M \alpha))$.
    These also commute with $\partial$ and pass to homology in the
    reduced equivariant theory.  A conjecture about $c_2,c_3,...$
    $s_2,s_3,...$ generating genus zero closed string operators
    and the algebraic form of this structure was proposed in [3] and is mentioned below in the
    summary. Also, compare [Chas] for the original questions
    motivating this work.

    \textit{Interplay between open and closed string states}: Let
    $\mathcal{C}$ denote the closed string states in $M$, a manifold of
    dimension $d$, and let $\vartheta$ denote any of the complexes of open
    string states.  Transversality yields an action of closed
    strings on open strings,

    \begin{center}
    $\mathcal{C} \otimes \vartheta \rightarrow \vartheta$  $\quad\quad\quad$     degree=($-d + 2$)
    \end{center}

     and a
    coaction of closed strings on open strings
    \begin{center}
     $\vartheta \rightarrow \mathcal{C} \otimes \vartheta$  $\quad\quad\quad$      degree=($-d + 2$)

    \end{center}

    In the coaction we let the open string hit itself at any two times
    and split the event into a closed string and an open string.  In
    the action we let a closed string combine with an open string to
    yield an open string.

    The action is a Lie action of the Lie algebra of closed strings by derivations at the
    transversal chain level.  Both the action and the coaction have a
    non trivial commutator with the boundary operator on chains.

    \S3\textit{Connection to work of (Atiyah-Segal),(Moore-Segal) and
    (Getzler and Segal)}:  Dialgebras satisfying the module or
    Frobenius compatibility give examples of 1+1 TQFT's in the positive
    boundary sense.  In the commutative case we associate the
    underlying vector space to a directed circle, its tensor products
    to a disjoint union of directed circles and to a connected $2D$
    oriented bordism between two non empty collections the morphism
    obtained by decomposing the bordism into pants and composing
    accordingly the algebra or coalgebra map.  The module or Frobenius
    compatibility is just what is required for the result to be
    independent of the choice of pants decomposition.

   $N.B.$ this description differs from the usual one because we don't have disks
   to close up either end of the bordism. One knows these discs
    at both ends would force the algebra to be finite dimensional and the algebra and
   coalgebra to be related by a non degenerate inner product.  We
   refer to these generalizations of the Atiyah-Segal concepts as the
   positive boundary version of TQFT (a name due to
   Ralph Cohen).

An exactly similar discussion with associative dialgebras
satisfying the module or Frobenius compatibility leads to a
positive boundary version of a relative TQFT using open intervals.
Now the algebra and coalgebra are associated to 1/2 pants (a disc
with $\partial$ divided into six intervals-three (1/2 seams)
alternating with two (1/2 cuffs) and one (1/2 waist)). Any planar
connected bordism between two nonempty collections of intervals
determines a mapping between inputs and outputs.

The structures we have found (including $\partial$ labels $L_a,
L_b,...$) for open strings using the composition $\wedge$ and
fixed time cutting $\vee_t$ satisfies this Frobenius compatibility
up to chain homotopy and we can apply it at the homology level in
the relative TQFT scheme just mentioned. This fits with the work
of Moore-Segal [M].

As we begin to look at the chain homotopy coproduct $\vee$ the
derivation or infinitesimal bialgebra compatibility appears.
According to [Gan] the derivation or infinitesimal bialgebra
compatibility is related to the notion of module or Frobenius
compatibility via Koszul duality (see appendix).

Now we are entering into a third stage-the proposal of Segal (and
independently Getzler) enriching the earlier notion of TQFT by
chain complexes and chain homotopies.

Recall the free loop space above gives on the ordinary (chain)
homology level a (strong homotopy) commutative associative product
and a cocommutative coassociative coproduct (cutting at a fixed
time) satisfying the module or Frobenius compatibility.  This
together with the associative Frobenius category above for open
strings fits with the model [M].  In that model ordinary and
equivariant levels are not distinguished.

We saw that passing to the equivariant setting the product and the
cutting at variable time gave a Lie bialgebra in the reduced
theory.  According to [Gan] Lie dialgebras with Drinfeld
compatibility are related to commutative dialgebras with Frobenius
compatibility by Koszul duality (see appendix).

\S4 \textit{Summary}: We have described the part of the
interpretation of open and closed string field theory in topology
associated to the basic product and coproduct (and in the
equivariant setting certain implied $n$-variable splitting and
collision operators as in [3]).  The coproduct discussion has two
levels involving a coproduct $\vee_t$ and an associated chain
homotopy coproduct $\vee$.

We found the open string product and the coproduct $\vee_t$
satisfied the module or Frobenius compatibility on the level of
homology. In a setting where $\vee_0$ and $\vee_1$ were zero or
even deformable to zero, $\vee$ emerges as or can be deformed to a
coproduct commuting with $\partial$ and thus a coproduct $\vee$ on
homology of one higher degree.  Then a new compatibility with the
product is observed- the derivation or infinitesimal bialgebra
compatibility (true transversally).

Similarly for the closed string one has to consider the free loop
space in both the ordinary and equivariant versions. For the open
string with diagonal boundary conditions the relevant ordinary
(chains) homology of the free loop space becomes a (strong
homotopy) commutative dialgebra with the module or Frobenius
compatibility. Passing to the equivariant theory required for the
closed string interpretation and reducing to kill $\vee_0$ and
$\vee_1$ which makes $\vee$ commute with $\partial$, the product
coproduct pair becomes a Lie dialgebra with the derivation or
Drinfeld compatibility (equals Lie bialgebra). According to [Gan]
the associative and commutative dialgebras with the module or
Frobenius compatibility are respectively Koszul dual to the
associative and Lie dialgebras with the derivation or Drinfeld
compatibility.  This suggests that one of the structures will
intervene in descriptions of strong homotopy versions (in the
sense of Stasheff) of the dual structure (see appendix).

One can go further as discussed in [3] and visualize conjecturally
all the above collision and splitting operations of the closed
string theory $c_2,c_3,...,s_2,s_3,... $ defining on homology a
structure Koszul dual to the positive boundary version of the
Frobenius manifold structure described in [Manin].

The above is only a partial interpretation. The full
interpretation of open closed string field theory in topology
involves full families of arbitrary cutting and reconnecting
operations of a string in an ambient space $M$.  For closed curves
some full families of these operators were labelled
combinatorially by decorated even valence ribbon graphs obtained
by collapsing chords in [3]. There is a serious compactness issue
for the full families discussed there for realizing these in
algebraic topology. The issue is a correct computation of the
boundary. The problem has a parallel with renormalization in
Feynman graphs. For the compactness algebraic topology issue one
needs to associate operators to families of geometric graphs where
various subgraphs are collapsing. When all the components of the
collapsing subgraphs are trees there is no real problem as
discussed in [3]. Similarly for Feynman graphs it is my
understanding that if there were only tree collapses there is no
problem of renormalization.

In both cases algebraic topology transversality normal bundle and
Feynman graphs the loops in collapsing subgraphs cause the
problems.

In [3] we had to deal with some simple cases of one loop subgraph
collapses to treat the identities defining the Lie bialgebra (in
particular Drinfeld compatibility).  This lead to the idea of
using the Fulton MacPherson compactification of configuration
spaces to complete the discussion.  There is a normal bundle issue
related to transversality which requires more analysis to treat
the general $FM$ stratum.  However for disjoint union of graphs
with at most one loop per component this normal bundle for
transversality can be easily described as in [3].

Now we expect a Riemann surface discussion to be sufficient to
complete the string field theory transversality construction. This
will complete the definition of the operations for this
topological interpretation of open closed string field theory. The
idea is that 1)  general cutting and reconnecting operation on
strings is isomorphic to the change in level that occurs when
passing through a critical level of a harmonic function on a
Riemann surface and 2) geometrical ideas due to Thurston and then
Penner [P] allow an analysis of the combinatorial
compactifications of spaces of Riemann surfaces in terms of ribbon
graphs.

Thus if the transversality cutting and reconnecting operations of
the string field theory interpretations are organized by ribbon
graphs, then the compactness and transversality normal bundle
issues discussed in [3] can be treated for open and closed
strings. This is work in progress.

\textit{Appendix}: (dialgebras and compatibilities) Let us call a
linear space $V$ with  two maps $V \otimes
V\overset{\wedge}\rightarrow V$ and $V \overset{\vee}\rightarrow V
\otimes V$ a dialgebra. Associative dialgebra means $\wedge$ is
associative and $\vee$ is coassociative.  Commutative dialgebra
means besides being associative $\wedge$ and $\vee$ are symmetric.
Lie dialgebra means both maps are skew symmetric and that jacobi
and cojacobi hold.

In all these cases $V$ and $V \otimes V$ have module structures
over $V$ and there are two kinds of compatibilities between
$\wedge$ and $\vee$ relative to these.  We get six kinds of
structures (five appear in this paper, see table below) which are
examples of definitions of algebras over dioperads [Gan]. These
are structures whose generators and relations are described
diagrammatically by trees.

The familiar example of a compatibility studied by Hopf that
$\vee$ is a map of algebras (associative or commutative case but
not Lie) can only be described by a non tree diagram.

The compatibilities we consider here are \newline derivation
compatibility {\quad} $\vee (a \cdot b)=(\vee a)\cdot b + a \cdot
\vee(b)$ and
\newline module compatibility {\qquad}{\quad} $  \vee (a \cdot b)=\vee(a)\cdot b= a
\cdot \vee(b)$

Table with names of compatibility and/or structure and/or
examples.
\newline
\hspace*{1.8in} \parbox[c]{1.5in}{Module\\compatibility}
\hspace*{1in}\parbox[c]{2.0in}{Derivation\\compatibility}\\

\parbox[t]{1in}{Associative\\dialgebra}
\hspace*{0.5in}\parbox[t]{2.0in}{Frobenius
compatibility\\$\Leftarrow$ Frobenius algebra=\\associative
algebra with\\non degenerate invariant\\inner product}
\hspace*{0.5in}\parbox[t]{2.0in}{infinitesimal
bialgebra\\compatibility=\\infinitesimal
bialgebra\\(see Aguilar)}\\

\parbox[t]{1in}{Commutative\\dialgebra}
\hspace*{.5in}\parbox[t]{2.0in}{Frobenius
compatibility\\$\Leftarrow$commutative Frobenius\\algebra}
\hspace*{.5in}\parbox[t]{2.8in}{commutative
cocommutative\\infinitesimal bialgebra}\\

\parbox[t]{1in}{Lie\\dialgebra}
\hspace*{.5in}\parbox[t]{2.0in}{Frobenius
compatibility\\$\Leftarrow$Lie algebra with\\non degenerate
invariant\\inner product} \hspace*{0.5in}
\parbox[c]{2.8in}{Drinfeld compatibility\\=Lie
bialgebra}\\

Where the $\cdot$ refers to the algebra structure or the module
structure (which means in the associative case $a\cdot(b\otimes
c)=(a\cdot b)\otimes c,(a\otimes b)\cdot c=a\otimes(b\cdot c)$
\newline and in the Lie case $a\cdot(b \otimes c)=-(b \otimes c)\cdot a=[a,b]\otimes c +
b\otimes[a,c]$ where $[x,y]=\wedge(x\otimes y$).)

In [Gan] Koszul dual pairs are defined and there it is proved that
upper left and upper right are Koszul dual pairs and that middle
left and lower right are Koszul dual pairs.  We suppose that the
lower left and middle right are also Koszul dual pairs.

We note in passing a remark about derivation or Drinfeld
compatibility and algebra or Hopf compatibility. A category of
"power series" Hopf algebras $\mathcal{U}$ was shown to be
equivalent to the category of Lie bialgebras $\mathcal{D}$ where
$\mathcal{D}\rightarrow\mathcal{U}$ was a formal quantization and
$\mathcal{U}\rightarrow\mathcal{D}$ was a semi classical limit
(Etingof-Kahzdan).

We emphasize these Koszul relations because in several important
situations a strong homotopy algebraic structure of one kind is
very naturally expressed by freely generated diagrams decorated
with tensors labeled by the Koszul dual structure.  In the above
discussion all the structures that are true transversally will
almost certainly lead to strong homotopy versions on the entire
space of states. So these might be expressed in this graphical
Koszul dual way.

\newpage
\textit{References}:

[1]K. Fukaya,Oh, Ohta, Ono "Lagrangin intersection Floer theory
-anomaly and obstruction" - (2000) See Fukaya website.

[2] M. Chas and D. Sullivan "String Topology" GT/ 9911159.
 Annals of Mathematics (to appear).

[3] M. Chas and D. Sullivan "Closed string operators in topology
leading to Lie bialgebras and higher string algebra" GT/ 0212358.
Abel Bicentennial Proceedings (to appear).

[G] E.Getzler "Operads and moduli spaces of genus zero Riemann
surfaces".  In:The Moduli Spaces of Curves ed. by R.Dijkgraaf, C.
Faber, G. van der Geer Progress in Math, vol.129 Birkhauser 1995,
199-230.

[M] Greg Moore "Some Comments on Branes, G-Flux, and K-theory"
Part II and references to Segal notes therein. International
Journal of Modern Physics A. arXiv:hep-th/0012007 v1. 1 Dec. 2000.

[Gan] Wee Liang Gan "Koszul duality for dioperads" preprint
University of Chicago 2002 QA/0201074.

[Goldman] William M.Goldman "Invariant functions on Lie group and
Hamiltonian flows of surface group representations" Invent. Math.
85(1986), no. 2, 263-302.

[Manin] Yuri Manin "Frobenius Manifolds, Quantum cohomology and
moduli spaces" AMS Colloquium Publications, Vol. 47.

[Cohen-Jones] "A homotopy theoretic realization of string
topology" math GT/0107187.

[Tradler] "The BV Algebra on Hochschild Cohomology Induced by
Infinity Inner Products", GT/0210150.

[S] James Stasheff "H-spaces from a homotopy point of view,"
Lecture Notes in Mathematics 161, Springer-Verlag, Berlin (1970),
ii-95.

[Chas] Moira Chas "Combinatorial Lie bialgebras of curves on
surfaces" to appear in Topology. Also arXiv GT/0105178.

[P] R. C. Penner "The decorated Teichmuller space of punctured
surface", Communications in mathematical physics 113 (1987)
299-339.

\quad
\begin{center}
CUNY Graduate Center, 365 Fifth Avenue, New York, NY  10016
\newline SUNY at Stony Brook, Stony Brook, NY  11794-3651
\newline  \textit{email: dsullivan@gc.cuny.edu}
\end{center}

    \end{document}